\documentclass[]{amsart}

\usepackage{amsmath}

\theoremstyle{plain}
\newtheorem{thm}{Theorem}[subsection]
\newtheorem{prop}[thm]{Proposition}
\newtheorem{lemma}[thm]{Lemma}

\theoremstyle{definition}

\theoremstyle{remark}
\newtheorem{rmk}[thm]{Remark}

\newcommand{\vol}{\ensuremath{\operatorname{vol}}}
\newcommand{\real}{\ensuremath{\operatorname{Re}}}
\newcommand{\imag}{\ensuremath{\operatorname{Im}}}
\newcommand{\spn}{\ensuremath{\operatorname{span}}}

\newcommand{\G}{\ensuremath{\operatorname{G_2}}}
\newcommand{\Gs}{\ensuremath{\operatorname{G_2}}{-structure}}
\newcommand{\SP}{\ensuremath{\operatorname{Spin (7)}}}
\newcommand{\SPs}{\ensuremath{\operatorname{Spin (7)}}{-structure}}

\newcommand{\ph}{\ensuremath{\varphi}}

\newcommand{\wtwop}{\ensuremath{\wedge^2_{+}}}
\newcommand{\wtwom}{\ensuremath{\wedge^2_{-}}}

\newcommand{\hk}{\ensuremath{\lrcorner}}

\newcommand{\tr}{\ensuremath{\operatorname{Tr}}}

\newcommand{\rest}[2]{\ensuremath{{#1} |_{{}_{#2}}}}
\newcommand{\be}{\ensuremath{\bar e}}
\newcommand{\ce}{\ensuremath{\check e}}
\newcommand{\bn}{\ensuremath{\bar \nu}}
\newcommand{\cn}{\ensuremath{\check \nu}}

\newcommand{\co}{\ensuremath{\check \omega}}
\newcommand{\nup}{\ensuremath{\nu^{\perp}}}
\newcommand{\oo}{\ensuremath{\mathbf 1}}
\newcommand{\oi}{\ensuremath{\mathbf i}}
\newcommand{\oj}{\ensuremath{\mathbf j}}
\newcommand{\ok}{\ensuremath{\mathbf k}}
\newcommand{\oee}{\ensuremath{\mathbf e}}
\newcommand{\oie}{\ensuremath{\mathbf i \mathbf e}}
\newcommand{\oje}{\ensuremath{\mathbf j \mathbf e}}
\newcommand{\oke}{\ensuremath{\mathbf k \mathbf e}}
\newcommand{\spp}{\ensuremath{{\, \slash \! \! \! \! \mathcal
S}_{\! \!{{}^+}}}}
\newcommand{\spm}{\ensuremath{{\, \slash \! \! \! \! \mathcal
S}_{\! \!{{}^-}}}}
\newcommand{\spi}{\ensuremath{{\, \slash \! \! \! \! \mathcal S}}}
\newcommand{\cliff}{\ensuremath{\! \! \cdot}}
\newcommand{\jm}{\ensuremath{\mathbf j_{\!{}_{\scriptscriptstyle
M}}}}

\numberwithin{equation}{section}
\numberwithin{table}{section}
\numberwithin{figure}{section}

\begin{document}

\title[Bundle Constructions of Calibrated Submanifolds in
$\mathbb R^7$ and $\mathbb R^8$]{Bundle Constructions of \\
Calibrated Submanifolds in $\mathbb R^7$ and $\mathbb R^8$}

\author{Marianty Ionel}
\address{Mathematics Department\\ McMaster University}
\email{ionelm@math.mcmaster.ca}
\author{Spiro Karigiannis}
\thanks{Spiro Karigiannis was partially supported by NSERC}
\address{Mathematics Department\\ McMaster University}
\email{spiro@math.mcmaster.ca}
\author{Maung Min-Oo}
\thanks{Maung Min-Oo was partially supported by NSERC}
\address{Mathematics Department\\ McMaster University}
\email{minoo@math.mcmaster.ca}
\subjclass{53, 58}
\date{\today}

\begin{abstract}
We construct calibrated submanifolds of $\mathbb R^7$
and $\mathbb R^8$ by viewing them as total spaces of
vector bundles and taking appropriate sub-bundles which are
naturally defined using certain surfaces in $\mathbb R^4$. We
construct examples of associative and coassociative submanifolds
of $\mathbb R^7$ and of Cayley submanifolds of $\mathbb R^8$. This
construction is a generalization of the Harvey-Lawson bundle
construction of special Lagrangian submanifolds of $\mathbb
C^{n}$.
\end{abstract}

\maketitle

\section{Introduction} \label{introsec}

The study of calibrated geometries was first initiated by Harvey
and Lawson 1982 in their seminal paper~\cite{HL}. Because they are
believed to play a crucial role in explaining the phenomenon of
mirror symmetry~\cite{SYZ}, they have recently received much
attention. There has been extensive research done on special
Lagrangian submanifolds of $\mathbb C^n$, most notably by Joyce
but see also~\cite{J2} and the many references contained therein.
Significantly less progress has been made in analyzing associative
and coassociative submanifolds of
$\mathbb R^7$ and Cayley submanifolds of $\mathbb R^8$, although
the recent papers~\cite{L1, L2} of Lotay presented some
constructions analogous to earlier special Lagrangian
constructions by Joyce. Needless to say, even less is known about
calibrated submanifolds in more general Calabi-Yau, \G, and \SP\
manifolds, even in the non-compact case, although the examples in
$\mathbb R^n$ serve as important local models, especially for
studying the possible singularities that can occur.

In their original paper~\cite{HL} Harvey and Lawson presented a
construction of special Lagrangian submanifolds in $\mathbb C^n$
using bundles. In this paper, motivated by their work, we describe
a similar bundle construction of associative and coassociative
submanifolds of $\mathbb R^7$ and Cayley submanifolds of $\mathbb
R^8$. The reader can consult~\cite{HL2, HL, J3} for background on
these exceptional calibrations.

The Harvey-Lawson contruction involves viewing $\mathbb C^n$ as a
vector bundle over $\mathbb R^n$, and taking an appropriate
sub-bundle of the restriction of this bundle to a submanifold $M^p
\subset \mathbb R^n$. In this case $\mathbb C^n = T^* (\mathbb
R^n)$ and the subbundle is the conormal bundle $N^* (M^p)$. They
find that the conormal bundle is special Lagrangian if and only in
$M^p$ is {\em austere} in $\mathbb R^n$, which is a condition
which is in general much stronger than minimal. Their construction
is reviewed in detail in Section~\ref{HLreviewsec}.

It is well known~\cite{BS} that if one views $\mathbb R^7$ as the
space of anti-self dual $2$-forms on $\mathbb R^4$, and $\mathbb R^8$
as the negative spinor bundle of $\mathbb R^4$, there are naturally
defined parallel \G\, and \SPs s on them, respectively. See~\cite{J1,
J3, K} for background on \G\ and \SPs s. We consider restricting these
bundles to a surface $M^2 \subset \mathbb R^4$, and then take
appropriate naturally defined sub-bundles of this restriction, the
total spaces of which are candidates for associative, coassociative,
and Cayley submanifolds. This is discussed in
Section~\ref{exceptionalsec}.

Since a calibrated submanifold is necessarily minimal, and since
the vector bundle directions have trivial second fundamental form,
the base manifold $M^2$ must be necessarily at least minimal in
$\mathbb R^4$. (Just as austere submanifolds are at least minimal
in the Harvey-Lawson construction.) In Theorem~\ref{coassthm} we
find that the naturally defined rank $2$ sub-bundle of
$\rest{\wedge^2_- (\mathbb R^4)}{M^2}$ is coassociative iff
the immersion of $M^2$ in $\mathbb R^4$ is a solution of
exactly one half of the {\em real isotropic minimal} surface
equation, sometimes also called superminimal. It is important that
not all real isotropic minimal surfaces will work. Perhaps
somewhat surprisingly, we find in Theorem~\ref{assthm} that the
naturally defined rank $1$ sub-bundle of $\rest{\wedge^2_- (\mathbb
R^4)}{M^2}$ is associative iff $M^2$ is just minimal in $\mathbb
R^4$, with no extra conditions. Similarly in
Theorem~\ref{cayleythm} we find two naturally defined rank $2$
sub-bundles of $\rest{\spm (\mathbb R^4)}{M^2}$ and each of them
is Cayley iff $M^2$ is again just minimal.

The associative construction produces interesting new examples,
while the coassociative construction actually produces examples
which live in a $\mathbb C^3$ subspace of $\mathbb R^7$ and are
complex submanifolds of $\mathbb C^3$. The Cayley construction
produces submanifolds of $\mathbb R^8$ which are either of the
form $\mathbb R \times L$ for an associative $3$-fold $L$ or are
non-trivial coassociative submanifolds of $\mathbb R^7$.

It is perhaps interesting that special Lagrangian and
coassociative submanifolds are harder to construct using these
methods, requiring a base manifold which is more than just
minimal. Special Lagrangian and coassociative submanifolds have a
very nice, unobstructed local deformation theory~\cite{M}, and the
local moduli space is intrinsic to the submanifold. On the other
hand, associative and Cayley submanifolds have a more complicated,
non-intrinsic and obstructed deformation theory, and yet the
bundle construction in these two cases is simpler, requiring only
minimality.

There exist examples of special holonomy metrics on non-compact
manifolds which are bundles over a compact base, for example the
Calabi-Yau metrics on $T^* (S^n)$, described in~\cite{St} and the
\G\ holonomy metrics on $\wedge^2_- (S^4)$ and $\wedge^2_- (\mathbb
C \mathbb P^2 )$ and the \SP\ holonomy metrics on $\spm (S^4)$,
described in~\cite{BS, GPP}. Similar constructions of calibrated
submanifolds can be done in these cases, and this is the
subject of a forthcoming paper~\cite{KM}.

We should remark that after this work was done, the authors
found a similar although different statement, without proof, in an
unpublished preprint by S.H. Wang~\cite{W}. His statement concerned
the non-compact \G\ and \SP\ manifolds first constructed by Bryant
and Salamon~\cite{BS} (we will deal with this case in~\cite{KM})
and his claim is that superminimal is the required condition for
all three constructions. We have proved that in the associative
and Cayley cases, just minimal is enough, while in the
coassociative case, we prove that only half of the superminimal
surfaces work.

{\em Acknowledgements.} The authors would like to thank Vestislav
Apostolov, Robert Bryant, and Mackenzie Wang for helpful
discussions.

\section{The Second Fundamental Form for Immersions $M^p \subset \mathbb
R^n$}
\label{computationssec}

In this section we set up notation for local computations for an
isometric immersion of a $p$-dimensional submanifold $M^p$ immersed in
$\mathbb R^n$.

Take $(x^1, x^2, \ldots, x^n)$ to be coordinates on $\mathbb R^n$, and
denote the immersion $M \subset \mathbb R^n$ by $x^i = x^i(u^1, u^2,
\ldots, u^p)$, for $1 \leq i \leq n$ where $(u^1, u^2, \ldots, u^p)$
are local coordinates on $M^p$. Consider a point $\mathbf u_0$ in $M$
with coordinates $(u^1_0, u^2_0, \ldots, u^p_0)$ and corresponding to
the point $\mathbf x_0 = \mathbf x (\mathbf u_0)$ in $X$ with
coordinates $(x^1_0, x^2_0, \ldots, x^n_0)$. Near $\mathbf x_0$ let
$e_1, e_2, \ldots , e_p$ be a local orthonormal frame of tangent
vector fields to $M$ and let $\nu_1, \nu_2, \ldots , \nu_q$ be a local
orthonormal frame of normal vector fields to $M$, where $q = n -
p$.

Let $\nabla$ denote the Levi-Civita connection on $\mathbb R^n$
and $( )^T$ and $( )^N$ denote the orthogonal projections onto the
tangent and normal bundles of $M$ in $\mathbb R^n$. By choosing an
orthonormal tangent frame and an orthonormal normal frame at the point
$\mathbf x_0$ and then parallel transporting via the induced tangent
and normal connections, we can assume that these local vector fields
have been chosen so that at the point $\mathbf x_0$,
\begin{equation} \label{zeroassumptionseq}
\rest{(\nabla_{e_i} e_j)}{\bf x_0}^T = 0 \qquad \text{and} \qquad
\rest{(\nabla_{e_i} \nu_j)}{\bf x_0}^N = 0 
\end{equation}

For $\nu$ any normal vector field, we can define the second
fundamental form $A^{\nu}$ as the linear operator
\begin{eqnarray*}
A^{\nu} & : & T(M) \to T(M) \\ A^{\nu} & : & w \mapsto A^{\nu}(w)
= {\left( \nabla_{w} \nu \right)}^T
\end{eqnarray*}
Here we are following the sign convention of Harvey and Lawson, which
differs from most definitions. The statements of all results in this
paper are independent of the choice of sign for the definition of
$A^{\nu}$. The important property of $A^{\nu}$ is that it is a
symmetric operator, and hence diagonalizable. This follows from
\begin{eqnarray*}
& & \langle e_i , A^{\nu}(e_j) \rangle = \langle e_i , {\left(
\nabla_{e_j} \nu \right)}^T \rangle = \langle e_i , 
\nabla_{e_j} \nu \rangle = - \langle \nabla_{e_j} e_i , \nu \rangle
\\ = & & - \langle \nabla_{e_i} e_j , \nu \rangle + \langle [ e_i,
e_j], \nu \rangle = - \langle \nabla_{e_i} e_j , \nu \rangle =
\langle e_j , A^{\nu}(e_i) \rangle
\end{eqnarray*}
where we have used $[e_i, e_j] = \nabla_{e_i} e_j - \nabla_{e_j}
e_i$ and the fact that $[e_i, e_j]$ is orthogonal to $\nu$ since
the bracket of two tangent vector fields on $M$ is again a tangent
vector field on $M$. We now adopt the notation
\begin{equation*}
A^{\nu}_{ij} = \langle A^{\nu} (e_i), e_j \rangle = A^{\nu}_{ji}
\end{equation*}
and more specifically $A^k_{ij} = A^{\nu_k}_{ij}$.

We also have the dual coframe of orthonormal cotangent vector fields
$e^1, e^2, \ldots, e^p$ and the orthonormal conormal vector fields
$\nu^1, \nu^2, \ldots, \nu^q$. These satisfy
\begin{equation} \label{dualformseq}
e^i(e_j) = \delta^i_j \qquad \nu^i(\nu_j) = \delta^i_j \qquad
e^i(\nu_j) = 0 \qquad \nu^i(e_j) = 0
\end{equation}
From~\eqref{dualformseq}, we have that $\left( \nabla_{e_i} e^j
\right) (e_k) = - e^j \left( \nabla_{e_i} e_k \right)$. From this, it
is very easy to check that under the hypotheses
of~\eqref{zeroassumptionseq}, we have the following expressions for
the covariant derivatives of the $e^i$'s and the $\nu^j$'s at the
point $\mathbf x_0$:
\begin{equation} \label{sffeq}
\nabla_{e_i} e^j = - \sum_{k=1}^q A^k_{ij} \nu^k \qquad \qquad
\nabla_{e_i} \nu^j = \sum_{k=1}^p A^j_{ik} e^k
\end{equation}

\section{The Harvey-Lawson Special Lagrangian Bundle
Construction}
\label{HLreviewsec}

In this section we review the bundle construction of Harvey and
Lawson~\cite{HL} of special Lagrangian submanifolds. The natural
ambient manifold in which to consider special Lagrangian
submanifolds is a Calabi-Yau manifold, which is in particular a
symplectic manifold. The simplest example of a symplectic manifold
is the {\em cotangent bundle} $T^* (\mathbb R^n)$ of $\mathbb
R^n$. This example is trivially Calabi-Yau, since
$T^* (\mathbb R^n) = \mathbb R^n \oplus \mathbb R^n = \mathbb C^n$.

On $\mathbb C^n = T^* (\mathbb R^n)$ we have a K\"ahler form $\omega =
\frac{i}{2} \sum dz^k \wedge d \bar z^k$ and a holomorphic $(n,0)$
volume form $\Omega = \real \Omega + i \imag \Omega = dz^1 \wedge
\ldots \wedge dz^n$. A real $n$-dimensional submanifold $L^n$ of
$\mathbb C^{n}$ is special Lagrangian with phase $e^{i \theta}$ (up to
a possible change of orientation) if the following two independent
conditions are satisfied:
\begin{equation*}
\rest{\omega}{L} = 0 \qquad \qquad \rest{\left( \imag e^{-i \theta}
\Omega \right)}{L} = 0
\end{equation*}
The first condition simply says that $L$ is
Lagrangian, which involves only the symplectic structure $\omega$
of $\mathbb C^n$. The {\em special} condition is given by the
second equation, which involves the Calabi-Yau metric structure.

Now it is a classical fact that if $M^p$ is a $p$-dimensional
submanifold of $\mathbb R^n$, then the {\em conormal bundle} $N^*
(M^p)$ is a {\em Lagrangian} submanifold of the symplectic
manifold $T^* (\mathbb R^n)$. (This will be shown below.) Motivated
by this, Harvey and Lawson found conditions on the immersion $M^p
\subset \mathbb R^n$ that makes $N^*(M^p)$ a special Lagrangian
submanifold of $T^* (\mathbb R^n)$, in terms of the second
fundamental form of the immersion. We reproduce their results
here, to motivate the constructions in Section~\ref{exceptionalsec}
and to fix our notation and conventions.

The canonical symplectic form $\omega$ on $T^* (\mathbb
R^n)$ is a $2$-form on the total space $T^* (\mathbb R^n) =
\mathbb R^n \oplus \mathbb R^n$. An orthonormal coframe for
$\mathbb R^n$ is given by $e^1, e^2, \ldots, e^n$. Hence an
arbitrary element of the cotangent bundle can be written as
\begin{equation*}
( \mathbf x, s_1 e^1 + s_2 e^2 + \ldots + s_n e^n )
\end{equation*}
where the $s_i$'s are coordinates on the cotangent space. An
orthonormal tangent frame for the total space is given by
\begin{eqnarray*}
(e_i , 0) \qquad i = 1, \ldots, n \qquad \text{ and } \qquad (0 ,
e^i) \qquad i = 1, \ldots, n
\end{eqnarray*}
For notational simplicity, we will denote $(e_i, 0)$ by $\be_i$ and
$(0, e^i)$ by $\ce^i$. The canonical symplectic form $\omega$ on
$T^* (\mathbb R^n)$ is then given by
\begin{equation*}
\omega = \sum_{k=1}^n \be^k \wedge \ce_k
\end{equation*}
where $\be^k$ is dual to $\be_k$ and $\ce^k$ is dual to $\ce_k$.  Let
$M^p \subset \mathbb R^n$. If we restrict the cotangent bundle $T^*
(\mathbb R^n)$ to $M^p$, we have
\begin{equation*}
\rest{T^* (\mathbb R^n)}{M} = T^* (M) \oplus N^* (M)
\end{equation*}
Since $M$ is $p$-dimensional, the total space of the conormal
bundle has dimension $p + (n-p) = n$. It therefore makes sense to
ask if $N^* (M)$ is Lagrangian.

We use the local coordinate notation described in
Section~\ref{computationssec}. An orthonormal coframe for
$\mathbb R^n$ is given by $e^1, e^2, \ldots, e^p, \nu^1, \nu^2,
\ldots, \nu^q$, where the $e_i$'s are tangent to $M^p$ and the
$\nu_i$'s are normal to $M^p$. Then $\omega$ takes the form
\begin{equation} \label{omegaeq}
\omega = \sum_{k=1}^p \be^k \wedge \ce_k + \sum_{l=1}^q \bn^l
\wedge \cn_l
\end{equation}
where as above $\bn_j = (\nu_j, 0)$ and $\cn^j = (0, \nu^j)$.

\begin{lemma} \label{slaglemma}
The conormal bundle $N^* (M)$ is a Lagrangian submanifold of $T^*
(\mathbb R^n)$.
\end{lemma}
\begin{proof}
We show that every tangent space to $N^* (M)$ is a Lagrangian
subspace of the corresponding tangent space to $T^* (\mathbb R^n)$.
In local coordinates the immersion $\Psi$ is given by
\begin{equation*}
\Psi : (u^1, u^2, \ldots, u^p, t_1, t_2, \ldots, t_q) \mapsto
(x^1(\mathbf u), \ldots, x^n(\mathbf u), t_1 \nu^1 + t_2 \nu^2 +
\ldots + t_q \nu^q )
\end{equation*}
Hence the tangent space at $(\mathbf x (\mathbf u_0), t_1, t_2,
\ldots, t_q)$ is spanned by the vectors
\begin{eqnarray*}
E_i & = & \Psi_* \left(\frac{\partial}{\partial u^i}\right) =
\left( e_i , \sum_{k=1}^q t_k \rest{\nabla_{e_i} (\nu^k)}{\mathbf
x_0} \right) \qquad i = 1, \ldots, p \\ F_j & = & \Psi_*
\left(\frac{\partial}{\partial t_j}\right) = (0, \nu^j) = \cn^j
\qquad j = 1, \ldots, q
\end{eqnarray*}
Using~\eqref{sffeq} we can write
\begin{equation*}
E_i = \left( e_i , \sum_{k=1}^q \sum_{l=1}^p t_k A^k_{il} e^l
\right) = \be_i + \sum_{l=1}^p A^{\nu}_{il} \ce^l
\end{equation*}
where we have defined $\nu = \sum_{k=1}^q t_k \nu_k$. To check that
the immersion is Lagrangian, we use \eqref{omegaeq} and compute
\begin{equation*}
\omega(F_i, F_j) = \omega( \cn^i, \cn^j ) = 0 \qquad
\forall i,j = 1, \ldots, q
\end{equation*}
and (dropping the summation sign over $k$ for clarity)
\begin{equation*}
\omega(F_i, E_j) = \omega( \cn^i, \be_j + A^{\nu}_{jk} \ce^k ) = 0
\qquad \forall i = 1, \ldots, q \quad j = 1, \ldots, p
\end{equation*}
Finally we have (again with the summations over $k$ and $l$ implied)
\begin{equation*}
\omega(E_i, E_j) = \omega( \be_i + A^{\nu}_{il} \ce^l, \be_j +
A^{\nu}_{jk} \ce^k ) =  A^{\nu}_{ij} - A^{\nu}_{ji} = 0
\end{equation*}
using the symmetry of $A^{\nu}$. Hence $\omega$ restricts to zero
on $N^* (M)$ and the conormal bundle is Lagrangian in $T^*
(\mathbb R^n)$.
\end{proof}

Since $T^* (\mathbb R^n)$ is Calabi-Yau, we can further
ask under what conditions the conormal bundle $N^* (M)$ is actually
{\em special Lagrangian}. A basis for the $(1,0)$ forms is given by
$\be^j + i \ce_j$ for $j = 1, \ldots, p$ and $\bn^k + i \cn_k$
for $k = 1, \ldots, q$. Thus the holomorphic $(n,0)$ form $\Omega$
can be written as
\begin{equation*}
\Omega = (\be^1 + i \ce_1) \wedge \ldots \wedge(\be^p + i \ce_p)
\wedge (\bn^1 + i \cn_1) \wedge \ldots \wedge (\bn^q + i \cn_q)
\end{equation*}

\begin{prop}[Harvey and Lawson, 1982~\cite{HL}, Theorem III.3.11] \label{slagprop}
The conormal bundle $N^* (M)$ is special Lagrangian in $T^* (\mathbb
R^n)$ with phase $i^q$ if and only if all the odd degree symmetric
polynomials in the eigenvalues of $A^{\nu}$ vanish for all normal
vector fields $\nu$ on $M$, where $A^{\nu}$ is the second fundamental
form for the immersion of $M$ in $\mathbb R^n$.
\end{prop}
\begin{rmk}
Such a submanifold is called {\em austere}. 
\end{rmk}
\begin{proof}
From Lemma~\ref{slaglemma} we had a basis for the tangent space
to the immersion of $N^* (M)$ at a point $(\mathbf x (\mathbf
u_0), t_1, t_2, \ldots, t_q)$ was given by
\begin{eqnarray*}
E_k = \be_k + \sum_{l=1}^p A^{\nu}_{kl} \ce^l & & k = 1, \ldots, p
\\ F_j = \cn^j & & j = 1, \ldots, q
\end{eqnarray*}
Without loss of generality we can assume that the tangent vector
fields were chosen to diagonalize $A^{\nu}$ at $\mathbf x_0$. That
is, $A^{\nu} (e_k) = \lambda_k e_k$ for $k = 1, \ldots, p$. We
compute easily that
\begin{equation*}
(\be^j + i \ce_j)(E_k) = \delta^j_k + i \lambda_k \delta^k_j \qquad \qquad
(\be^j + i \ce_j)(F_k) = 0
\end{equation*}
\begin{equation*}
(\bn^j + i \cn_j)(E_k) = 0 \qquad \qquad \qquad \qquad (\bn^j + i
\cn_j)(F_k) = \delta^j_k
\end{equation*}
and hence
\begin{equation*}
\Omega (E_1,\ldots, E_p, F_1, \ldots, F_q) = i^q (1 + i
\lambda_1)(1 + i \lambda_2) \cdots (1 + i \lambda_p)
\end{equation*}
If instead we consider the point $(\mathbf x (\mathbf u_0), c
t_1, c t_2, \ldots, c t_q)$ then the eigenvalues of $A^{c \nu}$
are $c \lambda_i$ and thus $\imag (i^{-q} \Omega)$ restricts to
zero on all these tangent spaces (for any $c$) if and only
if all the odd degree symmetric polynomials in the eigenvalues
vanish.
\end{proof}
\begin{rmk} \label{austerermk}
The first symmetric polynomial is the trace, so the submanifold
$M^p$ is necessarily minimal, as expected. If $p = 1,2$ this is
the only condition, but for $p \geq 3$ the austere condition is
much stronger than minimal.
\end{rmk}
\begin{rmk}
It is interesting to note that we cannot construct special
Lagrangian submanifolds in this way of arbitrary phase. The factor
of $i^{-q}$ means that the allowed phase (up to orientation)
depends on the codimension $q$ of the immersion.
\end{rmk}

\section{Bundle Constructions for Exceptional Calibrations}
\label{exceptionalsec}

Motivated by the results of Harvey and Lawson for constructing special
Lagrangian submanifolds using bundles, we look for a similar procedure
which will produce exceptional calibrated submanifolds: associative
and coassociative submanifolds of $\mathbb R^7$, and Cayley
submanifolds of $\mathbb R^8$.  The idea is as follows. There are
natural ways to view $\mathbb R^7$ and $\mathbb R^8$ as total spaces
of vector bundles over the base space $\mathbb R^4$, which are
compatible with the canonical \G\ and \SPs s on $\mathbb R^7$ and
$\mathbb R^8$. Specifically, the bundle of anti-self-dual $2$-forms
$\wtwom (\mathbb R^4) \cong \mathbb R^7$ has a natural \Gs, and the
negative spinor bundle $\spm (\mathbb R^4) \cong \mathbb R^8$ has a
natural \SPs. These structures are both parallel (torsion-free).

Now we let $M^p$ be a submanifold immersed in $\mathbb R^4$ and
consider the restriction of these bundles to $M^p$. For the right
choice of dimension $p$, this restriction breaks up naturally into
the direct sum of bundles, which can have the correct dimension
(as total spaces) to be candidates for calibrated submanifolds.
Then we can find conditions on the second fundamental form of the
immersion of $M$ in $\mathbb R^4$ for this to actually happen.
As discussed in Section~\ref{introsec}, we know that the
conditions on $M$ must include (at least) being a minimal
immersion.

\subsection{The space $\wtwom(\mathbb R^4)$ as a manifold with a
parallel \Gs}
\label{antiselfdualsec}

The space of anti-self-dual $2$-forms $\wtwom (\mathbb R^4)$
on $\mathbb R^4$ (which we will sometimes denote simply as
$\wtwom$) is naturally isomorphic to $\mathbb R^7$, with a natural
\Gs\ which we will now describe. (See~\cite{BS}, for example.) Let
$e^1, e^2, e^3, e^4$ be an oriented coframe of orthonormal covector
fields on $\mathbb R^4$. Then a basis of sections for $\wtwom$ is
given by
\begin{eqnarray*}
\omega^1 & = & e^1 \wedge e^2 - e^3 \wedge e^4 \\ \omega^2 & = &
e^1 \wedge e^3 - e^4 \wedge e^2 \\ \omega^3 & = & e^1 \wedge e^4 -
e^2 \wedge e^3
\end{eqnarray*}
The canonical $\G$ form $\ph$ on $\wedge^2_- (\mathbb
R^4)$ is a $3$-form on the total space $\wedge^2_- (\mathbb
R^4) = \mathbb R^4 \oplus \mathbb R^3$. An arbitrary element of
$\wedge^2_- (\mathbb R^4)$ can be written as
\begin{equation*}
\left( \mathbf x, t_1 \omega^1 + t_2 \omega^2 + t_3 \omega^3
\right) 
\end{equation*}
An orthonormal tangent frame for the total space is given by
\begin{eqnarray*}
(e_i , 0) \qquad i = 1, \ldots, 4 \qquad \text{ and } \qquad (0 ,
\omega^i) \qquad i = 1, \ldots, 3
\end{eqnarray*}
For notational simplicity, we will denote $(e_i, 0)$ by $\be_i$ and
$(0, \omega^i)$ by $\co^i$. The canonical $3$-form $\ph$ on
$\wedge^2_- (\mathbb R^4)$ is then given by
\begin{eqnarray} \label{g2formeq}
\ph & = & \co_1 \wedge \co_2 \wedge \co_3 + \co_1 \wedge ( \be^1
\wedge \be^2 - \be^3 \wedge \be^4 ) \\ & & {}+ \co_2 \wedge ( \be^1
\wedge \be^3 - \be^4 \wedge \be^2 ) + \co_3 \wedge ( \be^1 \wedge
\be^4 - \be^2 \wedge \be^3 )
\end{eqnarray}
where $\co_k$ is dual to $\co^k$ and $\be^k$ is dual to $\be_k$.
\begin{rmk} \label{signsrmk}
Alternatively, we could consider the bundle $\wtwop (\mathbb R^4)$ of
{\em self-dual} $2$-forms and obtain a \Gs\ on this space
where the three plus signs in~\eqref{g2formeq} become minus signs,
and the three minus signs become plus signs.
\end{rmk}

Let $M^2$ be a surface isometrically immersed in $\mathbb R^4$.
As in Section~\ref{computationssec}, we let $e_1, e_2$ be a local
orthonormal frame of tangent vector fields to $M$ and $\nu_1,
\nu_2$ be a local orthonormal frame of normal vector fields to $M$.
Then the dual covector fields $e^1, e^2$ and $\nu^1, \nu^2$ are
local coframes for the cotangent and conormal bundles. Locally we
can write that the anti-self-dual $2$-forms restrict to $M$ as
\begin{equation*}
\rest{\wtwom (\mathbb R^4)}{M} = \spn ( \omega^1, \omega^2, \omega^3)
\end{equation*}
where $\omega^1 = e^1 \wedge e^2 - \nu^1 \wedge \nu^2$, $\omega^2
=  e^1 \wedge \nu^1 - \nu^2 \wedge e^2$, and $\omega^3 =  e^1
\wedge \nu^2 - e^2 \wedge \nu^1$. Then $\omega^1$ is globally
defined on $M$ independent of the choice of orthonormal tangent
frames $e_1, e_2$ and normal frames $\nu_1$, $\nu_2$. Hence $\spn
(\omega^1)$ defines a rank $1$ bundle $E$ over $M^2$ and its
orthogonal complement (locally defined as $\spn (\omega^2,
\omega^3)$) defines a rank $2$ bundle $F$ over $M^2$.
\begin{equation*}
\rest{\wtwom (\mathbb R^4)}{M} = E \oplus F
\end{equation*}
The total spaces of $E$ and $F$ are $3$ and $4$-dimensional
submanifolds of $\mathbb R^7$ and hence candidates for
associative and coassociative submanifolds, respectively. Before
proceeding to check when this happens, we develop some formulas
that will be needed.

\begin{prop} \label{selfdualformulasprop}
Using the notation of Section~\ref{computationssec}, we have the
following expressions for the covariant derivatives of $\omega^1,
\omega^2, \omega^3$ in the $e_1, e_2$ directions at the point
$\mathbf x_0$.
\begin{eqnarray*}
\nabla_{e_i} \omega^1 & = & \left( A^2_{i1} - A^1_{i2} \right)
\omega^2 + \left( -A^1_{i1} - A^2_{i2} \right) \omega^3 \\
\nabla_{e_i} \omega^2 & = & \left( A^1_{i2} - A^2_{i1} \right)
\omega^1 \\ \nabla_{e_i} \omega^3 & = & \left( A^2_{i2} + A^1_{i1}
\right) \omega^1
\end{eqnarray*}
\end{prop}
\begin{proof}
We prove the second expression. We use~\eqref{sffeq} and compute:
\begin{eqnarray*}
\nabla_{e_i} \omega^2 & = & \left(\nabla_{e_i} e^1 \right)  \wedge
\nu^1 + e^1 \wedge \left(\nabla_{e_i} \nu^1 \right) -
\left(\nabla_{e_i} \nu^2 \right) \wedge e^2 - \nu^2 \wedge
\left(\nabla_{e_i} e^2 \right) \\ & = & \left( - A^1_{i1} \nu^1 -
A^2_{i1} \nu^2 \right) \wedge \nu^1 + e^1 \wedge \left(
A^1_{i1} e^1 + A^1_{i2} e^2 \right) \\ & & {}- \left( A^2_{i1} e^1
+ A^2_{i2} e^2 \right) \wedge e^2 - \nu^2 \wedge \left( -
A^1_{i2} \nu^1 - A^2_{i2} \nu^2 \right) \\ & = &
\left( A^1_{i2} - A^2_{i1} \right) \left( e^1 \wedge e^2 - \nu^1
\wedge \nu^2 \right)
\end{eqnarray*}
The other two are obtained similarly.
\end{proof}

\subsection{Coassociative Submanifolds of $\wtwom(\mathbb R^4)$}
\label{coasssec}

We are now ready to determine conditions on the immersion $M^2 \subset
\mathbb R^4$ so that the total space of the bundle $F$ over $M$ is a
coassociative submanifold. A $4$-manifold $L^4$ is coassociative (see
\cite{HL2} and \cite{HL} Section IV.1.B) if and only if
$\rest{\ph}{L^4} = 0$ where $\ph$ is the $3$-form defining the
\Gs.

In anticipation of our results, we need to first make some
definitions. A rank $2$ real vector bundle which is both oriented
and possesses a Riemannian metric on each fibre comes equipped
with a natural almost complex structure $J$ defined as follows. If
$v_1, v_2$ is an oriented orthonormal basis in a fixed fibre, we
define $J v_1 = v_2$ and $J v_2 = - v_1$. If we change the
orientation, we change $J$ to $-J$. In our setting, both
$T (M)$ and $N (M)$ are rank $2$ vector bundles with induced
Riemannian metrics coming from the isometric immersion of $M$ into
$\mathbb R^4$. Since $\mathbb R^4$ is taken to be oriented,
even if $M$ is not oriented, a choice of orientation on $T (M)$
induces an orientation on $N (M)$ and vice-versa.

\begin{thm} \label{coassthm}
The total space of the rank $2$ bundle $F$ over $M$ is a
coassociative submanifold of $\wtwom (\mathbb R^4)$ if and only if
the second fundamental form $A^{\nu}$ of the immersion $M \subset
\mathbb R^4$ satisfies
\begin{equation} \label{conditioneq}
A^{J \nu} = -J A^{\nu}
\end{equation}
for all normal vector fields $\nu$.
\end{thm}
\begin{rmk}
In this equation the $J$ on the left hand side corresponds to the
natural almost complex structure on $N(M)$ while the $J$ on the right
hand side corresponds to the natural almost complex structure on $T
(M)$.  Explicitly, $A^{J \nu} (w) = -J( A^{\nu} (w) )$ for all tangent
vectors $w$. Note that this condition is independent of the choice of
orientation of $T (M)$, since it determines the orientation on $N (M)$
and changing $J$ to $-J$ in both sides of this equation leaves it
invariant.
\end{rmk}
\begin{proof}
We show that every tangent space to $F$ is a coassociative
subspace of the corresponding tangent space to $\wtwom$. In
local coordinates the immersion $\Psi$ is given by
\begin{equation*}
\Psi : (u^1, u^2, t_2, t_3) \mapsto (x^1(u^1, u^2), x^2(u^1, u^2), t_2
\omega^2 + t_3 \omega^3 )
\end{equation*}
Hence the tangent space at $(\mathbf x (\mathbf u_0), t_2, t_3)$ is
spanned by the vectors
\begin{eqnarray*}
E_i & = & \Psi_* \left(\frac{\partial}{\partial u^i}\right) =
\left( e_i , t_2 \rest{\nabla_{e_i} (\omega^2)}{\mathbf x_0} + t_3
\rest{\nabla_{e_i} (\omega^3)}{\mathbf x_0} \right) \qquad i = 1, 2
\\ F_j & = & \Psi_* \left(\frac{\partial}{\partial t_j}\right) =
( 0, \omega^j) = \co^j \qquad j = 2, 3
\end{eqnarray*}
Using Proposition~\ref{selfdualformulasprop} we can write
\begin{eqnarray*}
E_1 & = & \be_1 + \left( t_2 \left( A^1_{12} - A^2_{11} \right) + t_3
\left( A^2_{12} + A^1_{11} \right) \right) \co^1 \\
E_2 & = & \be_2 + \left( t_2 \left( A^1_{22} - A^2_{12} \right) + t_3
\left( A^2_{22} + A^1_{12} \right) \right) \co^1
\end{eqnarray*}
If we now define the vectors $\nu = t_2 \nu_1 + t_3 \nu_2$ and $\nup =
- t_3 \nu_1 + t_2 \nu_2$, which are orthogonal normal vectors, then
the expressions for $E_1, E_2$ simplify to
\begin{eqnarray*}
E_1 & = & \be_1 + \left( A^{\nu}_{12} - A^{\nup}_{11} \right) \co^1 \\
E_2 & = & \be_2 + \left( A^{\nu}_{22} - A^{\nup}_{12} \right) \co^1
\end{eqnarray*}
Now since we have
\begin{eqnarray*}
\ph & = & \co_1 \wedge \co_2 \wedge \co_3 + \co_1 \wedge (\be^1 \wedge
\be^2 - \bn^1 \wedge \bn^2 ) \\ & & {}+ \co_2 \wedge (\be^1 \wedge
\bn^1 - \bn^2 \wedge \be^2 ) + \co_3 \wedge (\be^1 \wedge \bn^2 -
\be^2 \wedge \bn^1 )
\end{eqnarray*}
we can check when the immersion is coassociative by finding when
$\ph$ restricts to zero on each of these tangent spaces. It is
easy to compute that
\begin{equation*}
\ph (E_1, E_2, \cdot) = E_2 \hk E_1 \hk \ph = \co_1
+ \left( \cdots \right) \be^1 + \left( \cdots \right) \be^2
\end{equation*}
and hence since $F_j = \co^j$ we see that $\ph(E_1, E_2, F_2) =
\ph(E_1, E_2, F_3) = 0$ always. It remains to check when $\ph
(F_2, F_3, E_j) = 0$ for $j = 1,2$. Since $\ph(F_2, F_3, \cdot) =
\co_1$ these become the pair of conditions
\begin{equation*}
A^{\nu}_{12} - A^{\nup}_{11} = 0 \qquad \qquad A^{\nu}_{22} -
A^{\nup}_{12} = 0
\end{equation*}
for the tangent space at $(\mathbf x_0, t_2, t_3)$ to be
coassociative. We get two more equations that must be satisfied by
demanding that the tangent space at $(\mathbf x_0, -t_3, t_2)$
also be coassociative. This corresponds to changing $t_2 \mapsto
-t_3$ and $t_3 \mapsto t_2$ in the above equations, which is
equivalent to $\nu \mapsto \nup$ and $\nup \mapsto -\nu$. This
gives
\begin{equation*}
A^{\nup}_{12} + A^{\nu}_{11} = 0 \qquad \qquad A^{\nup}_{22} +
A^{\nu}_{12} = 0
\end{equation*}
Thus we see that at each point $\mathbf x (\mathbf u_0)$ on the
surface $M^2$, the matrix $A^{\nup}$ is determined by $A^{\nu}$ for
all normal vector fields $\nu$. We can combine the above four
equations in the following matrix equation:
\begin{equation} \label{coassconditionseq}
\begin{pmatrix} A^{\nup}_{11} & A^{\nup}_{12} \\ A^{\nup}_{12} &
A^{\nup}_{22} \end{pmatrix} = \begin{pmatrix} A^{\nu}_{12} &
A^{\nu}_{22} \\ -A^{\nu}_{11} & -A^{\nu}_{12} \end{pmatrix} =
\begin{pmatrix} 0 & 1 \\ -1 & 0 \end{pmatrix} \begin{pmatrix}
A^{\nu}_{11} & A^{\nu}_{12} \\ A^{\nu}_{12} & A^{\nu}_{22}
\end{pmatrix}
\end{equation}
which says $A^{\nup} = A^{J \nu} = -J A^{\nu}$ for $J =
\begin{pmatrix} 0 & -1 \\ 1 & 0 \end{pmatrix}$ which is the natural
almost complex structure described above. 
\end{proof}
\begin{rmk} \label{coassconditionsrmk2}
It is easy to check that if we assume that the second fundamental
forms $A^{\nu}$ and $A^{\nup}$ are {\em simultaneously
diagonalizable} and satisfy $A^{J \nu} = -J A^{\nu}$
then necessarily $A^{\nu} = A^{\nup} = 0$ and $M^2$ is {\em totally
geodesic} in $\mathbb R^4$, and hence is a plane. The constructed
coassociative submanifold is then a $4$-plane.
\end{rmk}

Note~\eqref{coassconditionseq} implies that $A^1_{11} + A^1_{22} =
A^2_{11} + A^2_{22} = 0$. Since $\nu$ and $\nup$ are a basis for the
normal space at every point, we see that $\tr (A) = 0$ and $M^2$ is
necessarily {\em minimal} in $\mathbb R^4$, as expected. However the
condition $A^{J \nu} = -J A^{\nu}$ is actually stronger than minimal,
just as the austere condition in Proposition~\ref{slagprop} was
stronger than minimal. These surfaces are well known, and are
sometimes called {\em superminimal}, although following the suggestion
of R.L. Bryant we prefer to call them {\em real isotropic minimal
surfaces}. We now give a more invariant description of these surfaces.

The second fundamental form $A$ can be viewed as a symmetric tensor
on $M^2$ with values in the normal bundle $N(M^2)$. That is,
\begin{equation*}
A = A_{11} e^1 \otimes e^1 + A_{12} e^1 \otimes e^2 + A_{21} e^2
\otimes e^1 + A_{22} e^2 \otimes e^2
\end{equation*}
where $A_{ij} = A_{ji} = A^1_{ij} \nu_1 + A^2_{ij} \nu_2$ and $e_1,
e_2$ and $\nu_1, \nu_2$ are oriented orthonormal tangent and normal
frames for $M^2$, respectively. We have a natural almost complex
structure $J$ on $T(M)$ given by $J e_1 = e_2$, and $J e_2 = -
e_1$. Therefore we can consider the $1$-forms $\beta = e^1 + i e^2$
and $\bar \beta = e^1 - i e^2$, which are of type $(1,0)$ and $(0,1)$
respectively.  We can rewrite $A$ in this basis as follows:
\begin{equation*}
A = \frac{\mathbf W}{4} \beta \otimes \beta + \frac{\mathbf H}{4}
(\beta \otimes \bar \beta + \bar \beta \otimes \beta) +
\frac{\overline{\mathbf W}}{4} \bar \beta \otimes \bar \beta
\end{equation*}
where $\mathbf H = A^{1}_{11} \nu_1 + A^{2}_{11} \nu_2 +
A^{1}_{22} \nu_1 + A^{2}_{22} \nu_2$ is the mean curvature
vector of the immersion and $\mathbf W$ is the complex
valued normal vector
\begin{equation*}
\mathbf W = \mathbf W_1 + i \mathbf W_2 = \left( A^{1}_{11} \nu_1 +
A^{2}_{11} \nu_2 - A^{1}_{22} \nu_1 - A^{2}_{22} \nu_2 \right) + i
\left( - 2 A^{1}_{12} \nu_1 - 2 A^{2}_{12} \nu_2 \right)
\end{equation*}
If $M$ is minimal in $\mathbb R^4$, then $\mathbf H = 0$ and the
$(1,1)$ term in $A$ vanishes. The $(2,0)$ and $(0,2)$ terms are
conjugates of each other. We use the real inner product on
$\mathbb R^4$ to define the {\em complex quartic form}
$Q$ of the minimal surface $M$ to be
\begin{equation*}
Q = \mathbf W \cdot \mathbf W = \left( \mathbf W_1 \cdot \mathbf
W_1 - \mathbf W_2 \cdot \mathbf W_2 \right) + 2 i \left( \mathbf
W_1 \cdot \mathbf W_2 \right)
\end{equation*}
The minimal surface $M$ is called {\em real isotropic} if $Q = 0$.
It is easy to check that real isotropic is equivalent
(using the fact that $M$ is already minimal) to the equations
\begin{equation*}
A^1_{11} = \pm A^2_{12} \qquad A^1_{22} = \mp A^2_{12} \qquad
A^2_{11} = \mp A^1_{12} \qquad A^2_{22} = \pm A^1_{12}
\end{equation*}
which can be written concisely as
\begin{equation} \label{supermineq}
A^{J \nu} = \pm J A^{\nu}
\end{equation}
for any normal vector field $\nu$. The
condition~\eqref{conditioneq} we obtained above was this equation
with {\em only} the minus sign. (And it appears with only the plus
sign if we are considering $\wedge^2_+ (\mathbb R^4)$.) Hence only
half of the real isotropic minimal surfaces in $\mathbb R^4$ can be
used in the construction. Real isotropic surfaces have been
extensively studied by many, and the interested reader can refer
to~\cite{Br3, ES, Sa2} and the references contained therein for
more details.

Suppose a surface $M^2 \subset \mathbb R^4$ satisfies $A^{J_2 \nu}
= -J_1 A^{\nu}$, where $J_1$ and $J_2$ are the natural almost
complex structures on the tangent and normal spaces,
respectively. (These were both referred to as $J$ above but now
we distinguish them explicitly for clarity.) We can define an
almost complex structure $\widetilde J$ on the rank $4$ vector
bundle $\rest{T^* (\mathbb R^4)}{M}$ over $M$ as follows: 
\begin{equation*}
\widetilde J = \begin{pmatrix} J_1 & 0 \\ 0 & - J_2 \end{pmatrix}
\end{equation*} 
acting diagonally on the tangent and normal spaces. In this
notation, the condition~\eqref{conditioneq} becomes $A^{\widetilde
J \nu} = \widetilde J A^{\nu}$. This is equivalent to
\begin{equation} \label{parallel1eq}
( \overline{ \nabla}_X (\widetilde J \nu) )^T =
\widetilde J ( \overline{ \nabla}_X \nu )^T
\end{equation}
where $\overline \nabla$ is the Levi-Civita connection on $\mathbb
R^4$, $X$ is a tangent vector field to $M$, and $\nu$ is a normal
vector field to $M$.

\begin{prop} \label{coassdegenerateprop}
If~\eqref{conditioneq} holds, then the almost complex structure
$\widetilde J$ defined above satisfies
\begin{equation*}
\overline \nabla_X \widetilde J = 0
\end{equation*}
for all tangent vector fields $X$ to $M$.
\end{prop}
\begin{proof}
Let $X$ and $Y$ be tangent vector fields to $M$. Using the fact
that $\widetilde J$ is orthogonal and also preserves the tangent
and normal spaces, we can use~\eqref{parallel1eq} to compute
\begin{eqnarray*}
\langle ( \overline{ \nabla}_X (\widetilde J \nu) )^T, Y \rangle &
= & \langle \widetilde J ( \overline{ \nabla}_X \nu )^T, Y \rangle
\\ - \langle \widetilde J \nu, \overline \nabla_X Y \rangle & = &
- \langle \overline \nabla_X \nu , \widetilde J Y \rangle \\ 
\langle \nu , \widetilde J ( \overline \nabla_X Y )^N \rangle & =
& \langle \nu, ( \overline \nabla_X \widetilde J Y )^N \rangle
\end{eqnarray*} 
which holds for all normal vector fields $\nu$, and hence
\begin{equation} \label{parallel2eq}
( \overline{ \nabla}_X (\widetilde J Y) )^N =
\widetilde J ( \overline{ \nabla}_X Y )^N
\end{equation}
Let $\nabla$ denote the Levi-Civita connection on $M^2$
from the induced metric, we have
\begin{eqnarray*}
\overline \nabla_X (\widetilde J) Y & = & \overline \nabla_X (
\widetilde J Y ) - \widetilde J (\overline \nabla_X Y ) \\ & = &
\nabla_X (\widetilde J Y) + ( \overline{ \nabla}_X
(\widetilde J Y) )^N - \widetilde J ( \nabla_X Y +
\left( \overline \nabla_X Y \right)^N ) \\ & = & \nabla_X ( J_1 Y
) - J_1 (\nabla_X Y) = \nabla_X (J_1) Y = 0 
\end{eqnarray*}
where we have used~\eqref{parallel2eq} in the third line and the
last equality is due to the fact that any almost complex structure
on a rank $2$ bundle is necessarily parallel. In the same
way~\eqref{parallel1eq} can be used to show
\begin{equation*}
\overline \nabla_X (\widetilde J) \nu = \nabla_X (J_2) \nu = 0 
\end{equation*}
and the result now follows.
\end{proof}
Unfortunately, Proposition~\ref{coassdegenerateprop} means that
all the coassociative submanifolds of $\mathbb R^7$ thus
constructed are everywhere orthogonal to a parallel direction,
given by $\omega_1 = e^1 \wedge e^2 - \nu^1 \wedge \nu^2$, and
actually live in an $\mathbb R^6$ subspace of $\mathbb R^7$. A
coassociative submanifold of $\mathbb R^7$ which misses one
direction is actually a complex dimension $2$ complex submanifold
of $\mathbb C^3 = \mathbb R^6$ (up to a possible change of
orientation). It is interesting to note, however, that precisely
which $\mathbb C^3$ sitting in $\mathbb R^7$ contains this complex
submanifold depends on the immersion of the surface $M^2$ in
$\mathbb R^4$.
\begin{rmk}
In Section~\ref{cayleysec}, during our search for Cayley
submanifolds of $\mathbb R^8$, we will obtain non-trivial
coassociative submanifolds of $\mathbb R^7$ which are not contained
in a strictly smaller subspace.
\end{rmk}
\begin{rmk}
On more general non-compact manifolds with holonomy \G\, such as
$\wedge^2_- (S^4)$ and $\wedge^2_- (\mathbb C \mathbb P^2)$
(see~\cite{BS, GPP}), this construction will produce more
interesting coassociative submanifolds. This is discussed
in~\cite{KM}.
\end{rmk}

\subsection{Associative Submanifolds of $\wtwom(\mathbb R^4)$}
\label{asssec}

Similarly we can determine conditions on the immersion $M^2 \subset
\mathbb R^4$ so that the total space of the bundle $E$ over $M$ is an
associative submanifold. A $3$-manifold $L^3$ is associative (see
\cite{HL2} and \cite{HL} Section IV.1.A) if and only if its tangent
space at every point $x$ is an associative subspace of $T_x (\wtwom
(\mathbb R^4) ) \cong \mathbb R^7$. Here we identify $\mathbb R^7
\cong \imag \mathbb O$, the imaginary octonions.

\begin{thm} \label{assthm}
The total space of the rank $1$ bundle $E$ over $M$ is an
associative submanifold of $\wtwom (\mathbb R^4)$ if and only if
the immersion $M \subset \mathbb R^4$ is {\em minimal}.
\end{thm}
\begin{proof}
We show every tangent space to $E$ is an associative
subspace of the corresponding tangent space to $\wtwom (\mathbb
R^4)$. In local coordinates the immersion $\Psi$ is
\begin{equation*}
\Psi : (u^1, u^2, t_1) \mapsto (x^1(u^1, u^2), x^2(u^1, u^2), t_1
\omega^1)
\end{equation*}
Hence the tangent space at $(\mathbf x (\mathbf u_0), t_1)$ is
spanned by the vectors
\begin{eqnarray*}
E_i & = & \Psi_* \left(\frac{\partial}{\partial u^i}\right) =
\left( e_i, t_1 \rest{\nabla_{e_i} (\omega^1)}{\mathbf x_0} \right)
\qquad i = 1, 2 \\ F_1 & = & \Psi_* \left(\frac{\partial}{\partial
t_1}\right) = (0, \omega^1) = \co^1
\end{eqnarray*}
From Proposition~\ref{selfdualformulasprop} we have
\begin{eqnarray*}
E_1 & = & \be_1 + t_1 \left( (A^2_{11} - A^1_{12}) \co^2 + (-A^1_{11}
- A^2_{12} ) \co^3 \right) \\ E_2 & = & \be_2 + t_1 \left( (A^2_{12} -
A^1_{22}) \co^2 + (-A^1_{12} - A^2_{22} ) \co^3 \right)
\end{eqnarray*}
To check that the tangent space at $(\mathbf x_0, t_1)$ is
associative, we need to verify that the associator $[E_1, E_2, F_1] =
(E_1 E_2)F_1 - E_1(E_2 F_1)$ vanishes. Without loss of generality, at
a point we can take the following explicit identification $T_x (\wtwom
(\mathbb R^4)) \cong \imag \mathbb O$:
\begin{equation*}
\begin{pmatrix} \co^1 & \co^2 & \co^3 & \be_1 & \be_2 &
\bn_1 & \bn_2 \\ \updownarrow & \updownarrow & \updownarrow &
\updownarrow & \updownarrow & \updownarrow & \updownarrow \\ \oi & \oj
& \ok & \oee & \oie & \oje & \oke \end{pmatrix}
\end{equation*}
and hence
\begin{eqnarray*}
E_1 & = & \oee + t_1 \left( (A^2_{11} - A^1_{12}) \oj + (-A^1_{11}
- A^2_{12} ) \ok \right) \\ E_2 & = & \oie + t_1 \left( (A^2_{12} -
A^1_{22}) \oj + (-A^1_{12} - A^2_{22} ) \ok \right) \\ F_1 & = &
\oi
\end{eqnarray*}
Now we can compute the associator (see the octonion multiplication
table in Appendix~\ref{octtablesec}), with the result being
\begin{eqnarray*}
[E_1,E_2,F_1] & = & (E_1 E_2) F_1 - E_1 (E_2 F_1) \\ & = & (-2
A^2_{11} - 2 A^2_{22}) \oje + (2 A^1_{11} + 2 A^1_{22}) \oke 
\end{eqnarray*}
which vanishes if and only if $\tr A^{\nu_1} = \tr A^{\nu_2} = 0$.
\end{proof}

\subsection{The space $\spm (\mathbb R^4)$ as a manifold with a parallel
\SPs}
\label{spinorssec}

The simplest \SPs\ on the total space of a bundle is the
negative spinor bundle $\spm (\mathbb R^4)$ of $\mathbb R^4$. We
will now explain how to see this. Over each point $x \in
\mathbb R^4$, the fibre of spinors over $x$ is isomorphic to two
copies of the quaternions $\spp \oplus \spm = \mathbb H
\oplus \mathbb H$. The one-forms (covectors) at $x$ are a
subset of the Clifford algebra over $x$, and hence act on the
spinor space. A good reference for spin representations is the
book of Harvey~\cite{H}. If $e^1, e^2, e^3, e^4$ is an orthonormal
basis of $1$-forms at $x$, then the Clifford algebra relations are
\begin{equation*}
e^i \cliff e^j + e^j \cliff e^i = -2 \delta^{i j}
\end{equation*}
where the $\cdot$ denotes the Clifford product. Clifford
multiplication by $1$-forms interchanges the two spaces
$\spi_{\pm}$. We identify the spinor space with the octonions, $\spp
\oplus \spm \cong \mathbb H \oee \oplus \mathbb H \cong \mathbb
O$. Octonionic multiplication by elements of $\mathbb H \oee$
interchanges $\mathbb H \oee$ and $\mathbb H$ (see
Appendix~\ref{octtablesec}). Also, we have the following identities
for octonionic multiplication (see \cite{HL} Appendix IV.A):
\begin{eqnarray*}
a (a x) & = & a^2 x \\ a_1( \bar {a_2} x) & = & - \bar {a_2} (a_1
x) \qquad \text{for } a_1, a_2 \text{ orthogonal}
\end{eqnarray*}
If we take $a_i \in \mathbb H \oee$, then $\bar {a_i} = -a_i$ and
hence if $e^1, e^2, e^3, e^4$ is an orthonormal basis of $\mathbb H
\oee$, these relations become
\begin{equation*}
e^i (e^j x) + e^j (e^i x) = - 2 \delta^{ij} x
\end{equation*}
Thus we see we obtain the spin representation at each point from
octonionic multiplication by identifying $\spp \oplus \spm \cong
\mathbb H \oee \oplus \mathbb H$ and the $1$-forms with $\mathbb H
\oee$. We will only require this representation for Clifford products
of $1$-forms and it will be written
\begin{eqnarray*}
\gamma & : & T^* \to \operatorname{End}(\spp \oplus \spm) \\
\gamma(\alpha) (s) & = & \alpha s 
\end{eqnarray*}
where $\alpha$ is a $1$-form, $s \in \spp \oplus \spm$ and the
product $\alpha s$ is octonionic multiplication. Note that since
$\mathbb O$ is {\em not associative}, we have to be careful when
composing two elements of this representation:
\begin{equation*}
\left(\gamma(\alpha_1) \gamma(\alpha_2)\right) (s) =
\gamma(\alpha_1) \left( \gamma(\alpha_2)(s) \right) = 
\gamma(\alpha_1) (\alpha_2 s) = \alpha_1 ( \alpha_2 s)
\end{equation*}
which in general is {\em not} the same as $(\alpha_1 \alpha_2) s$.

Now a manifold has a \SPs\ if at every point its tangent space can
be naturally identified with $\mathbb O$. With the identifications
we have made, the total space of $\spm (\mathbb R^4)$ has a tangent
space (at a point) isomorphic to $T(\mathbb R^4) \oplus \spm \cong
T^*(\mathbb R^4) \oplus \spm \cong \mathbb H \oee \oplus \mathbb H
\cong \mathbb O$.

Proceeding as before, we now isometrically immerse a
submanifold $M^p$ in $\mathbb R^4$ so that the restriction
$\rest{\spm (\mathbb R^4)}{M^p}$ splits naturally into pieces, and
hope to obtain Cayley submanifolds in this way. Once again, the
only natural choice occurs when $p=2$, the case of a surface. If
we let
$e^1, e^2$ be a local orthonormal coframe for $M^2$, and $\nu^1,
\nu^2$ a local orthonormal basis for the conormal bundle, then we
can consider the operations on the fibre $\spm$ of Clifford
multiplication with $\gamma(e^1) \gamma(e^2)$ or $\gamma(\nu^1)
\gamma(\nu^2)$. Two remarks are in order. First, since
multiplication by $\gamma(\alpha)$ interchanges $\spp$ and $\spm$,
we need to consider the composition of two such multiplications to
stay in $\spm$. Second, up to a sign (corresponding to a
choice of orientation for $M^2$) these operators are independent of
the choice of $e^1, e^2$ or $\nu^1, \nu^2$ since, for example
$\gamma(e^1)\gamma(e^2) = \gamma(e^1 \cliff e^2) = \gamma(e^1
\wedge e^2)$ because $e^1$ and $e^2$ are orthonormal.

The spinor space $\spm$ can be given the structure of a complex
$2$-dimensional vector space in many ways. One can check that if
$a, b, p, q \in \mathbb H$, then
\begin{equation*}
(a \oee) ((b \oee) (p q)) = p ( (a \oee) ((b \oee) q) )
\end{equation*}
That is, left multiplication by ordinary quaternions $\mathbb H$
commutes with the composition of {\em two} left multiplications by
elements of $\mathbb H \oee$. Now left multiplication by a unit
imaginary quaternion is a complex structure on $\mathbb H$, so
$\spm$ has an $S^2$ family of complex structures with respect to
operators of the form $\gamma (a\oee) \gamma(b\oee) : \spm \to
\spm$. Since we have a surface $M^2$ immersed in $\mathbb R^4$,
this determines a canonical complex structure $\jm$ on $\spm$
as follws. If $e^1 = a \oee$ and $e^2 = b \oee$ are an orthonormal
basis of tangent vectors to $M$, then $\jm$ is defined by
\begin{equation*}
\jm = e^1 e^2 = (a \oee) (b \oee) = - \bar b a
\end{equation*}
It is easy to check that $\jm$ is purely imaginary, and of
unit length, so $\jm^2 = -1$.  Alternatively, if we had used an
orthonormal basis of the normal space $\nu^1$ and $\nu^2$ and
multiplied them together as elements of $\mathbb H \oee$, we would
have obtained $-\jm$. Either choice will produce the same results
below.

\begin{lemma} \label{eigenspaceslemma}
The operator $r_T = \gamma(e^1)\gamma(e^2)$ satisfies $r_T^2 = -1$
and hence decomposes the space $\spm$ into two
$2$-dimensional eigenspaces $V_{\pm \jm}$ of eigenvalues $\pm \jm$.
Further, the operator $r_N = \gamma(\nu^1)\gamma(\nu^2)$ is equal
to $r_T$.
\end{lemma}
\begin{proof}
We compute $(r_T)^2 = \gamma(e^1 \cliff e^2
\cliff e^1 \cliff e^2) = - \gamma(1) = -1$ using the fact that $e^1
\cliff e^2 = - e^2 \cliff e^1$ and $e^i \cliff e^i = -1$. The
eigenspace decomposition now follows. Also,
$\gamma(e^1)\gamma(e^2)\gamma(\nu^1)\gamma(\nu^2) = \gamma(e^1
\cliff e^2 \cliff \nu^1 \cliff \nu^2) = \gamma(\vol)$ where $\vol$
is the volume form, and the spinor spaces $\spi_{\pm}$ are defined
as $\pm 1$ eigenspaces of Clifford multiplication with $\vol$:
$\vol \, \cliff \spi_{\pm} = \pm \spi_{\pm}$. Thus
$r_T r_N$ is minus the identity on $\spm$ and since $r_T$ and
$r_N$ commute (and hence are simultaneously diagonalizable), it is
easy to see that we must have $r_T = r_N$.
\end{proof}
We will henceforth denote $r_T = r$. Note that since we are only
interested in the eigenspaces $V_{\pm \jm}$, it does not matter
which orientation we choose for $M^2$. In fact, we can identify
these eigenspaces exactly. The octonion multiplication rules show
\begin{equation*}
(a\oee) ( (b \oee) q) = - \overline{ (b \bar q) } a = -
q \bar b a = q \jm
\end{equation*}
and so the operator $r$ is exactly right multiplication by $\jm$.
Thus the $+ \jm$ eigenspace of $r$ is $\spn { \{\oo, \jm \} }$ and
the $- \jm$ eigenspace is the orthogonal complement of this.

\subsection{Cayley Submanifolds of $\spm(\mathbb R^4)$}
\label{cayleysec}

We have described the natural splitting
\begin{equation*}
\rest{\spm (\mathbb R^4)}{M^2} = V_{+ \jm} \oplus V_{- \jm}
\end{equation*}
into two rank $2$ bundles over the base surface $M^2$. The total
space of either of these bundles is $4$-dimensional and is a
candidate for being a Cayley submanifold.

\begin{thm} \label{cayleythm}
The total space of either rank $2$ bundle $V_{\pm \jm}$ over $M$ is
a Cayley submanifold of $\spm (\mathbb R^4)$ if and only the
immersion $M \subset \mathbb R^4$ is {\em minimal}.
\end{thm}
\begin{proof}
We show every tangent space to the total space of $V_{+\jm}$ is a
Cayley  subspace of the corresponding tangent space to $\spm
(\mathbb R^4)$. The proof for $V_{-\jm}$ is identical. In local
coordinates the immersion
$\Psi$ is
\begin{equation*}
\Psi : (u^1, u^2, t_1, t_2) \mapsto (x^1(u^1, u^2), x^2(u^1, u^2), t_1
q_1(u^1,u^2) + t_2 q_2(u^1,u^2))
\end{equation*}
where $q_1$ and $q_2$ are an orthonormal basis of $V_{+\jm}$ and
hence satisfy $r q_k = \jm q_k$. The tangent space at $(\mathbf x
(\mathbf u_0), t_1,t_2)$ is spanned by the vectors
\begin{eqnarray*}
E_k & = & \Psi_* \left(\frac{\partial}{\partial u^k}\right) = e_k +
\rest{\nabla_{e_k} (t_1 q_1 + t_2 q_2)}{\mathbf x_0} \qquad k = 1,
2 \\ F_k & = & \Psi_* \left(\frac{\partial}{\partial t_k}\right) =
q_k \qquad k = 1, 2
\end{eqnarray*}
We now derive an expression for $\rest{\nabla_{e_k} q_j}{\mathbf
x_0}$. To simplify notation we will use a dot to denote
$\rest{\nabla_{e_k}}{\mathbf x_0}$. Since $r^2 = -1$, we can
differentiate to obtain
\begin{equation*}
r \dot r + \dot r r = 0
\end{equation*}
Hence since $r$ and $\dot r$ anti-commute, $r(\dot r q_j) = - \dot
r (r q_j) = - \jm \dot r q_j$ and thus $\dot r q_j \in V_{-\jm}$.
Now differentiating the equation $r q_j = \jm q_j$, we have
\begin{eqnarray*}
\dot r q_j + r \dot q_j & = & \jm \dot q_j \\ (r - \jm) \dot q_j &
= & - \dot r q _j
\end{eqnarray*}
The right hand side is in $V_{-\jm}$, and on this space $r = -\jm$,
so $r - \jm = -2\jm$ on $V_{-\jm}$ and we have
\begin{equation*}
(r - \jm)^{-1} (r - \jm) \dot q_j = \dot q_j = \frac{-1}{2}{(-\jm)}
(- \dot r q_j) = -\frac{\jm}{2} \dot r q_j
\end{equation*}
Explicitly, at the point $\mathbf x_0$, we have
\begin{equation*}
\nabla_{e_k} q_j = -\frac{\jm}{2} \left(
\gamma( \nabla_{e_k} e^1)\gamma(e^2) + \gamma(e^1)\gamma(
\nabla_{e_k} e^2) \right) q_j
\end{equation*}
From~\eqref{sffeq} this can be written as
\begin{equation*}
\nabla_{e_1} q_j = \frac{\jm}{2} \left( a_{11}
\gamma(\nu^1)\gamma(e^2) + b_{11}\gamma(\nu^2)\gamma(e^2) +
a_{12}\gamma(e^1)\gamma(\nu^1) + b_{12}\gamma(e^1)\gamma(\nu^2)
\right) q_j
\end{equation*}
\begin{equation*}
\nabla_{e_2} q_j = \frac{\jm}{2} \left( a_{12}
\gamma(\nu^1)\gamma(e^2) + b_{12}\gamma(\nu^2)\gamma(e^2) +
a_{22}\gamma(e^1)\gamma(\nu^1) + b_{22}\gamma(e^1)\gamma(\nu^2)
\right) q_j
\end{equation*}
where we have used the notation $a_{ij} = \langle e_i, A^{\nu_1} (e_j)
\rangle$ and $b_{ij} = \langle e_i, A^{\nu_2}(e_j) \rangle$.  Note
that the operators $\gamma(e^i)\gamma(\nu^j)$ all anti-commute with $r
= \gamma(e^1)\gamma(e^2)$ and hence map $V_{+\jm} \to
V_{-\jm}$. Therefore $\nabla_{e_k} q_j \in V_{-\jm}$.  To check that
the tangent space at $(\mathbf x_0, t_1, t_2)$ is Cayley, we need to
verify that the purely imaginary $4$-fold octonion product $\imag (E_1
\times E_2 \times F_1 \times F_2)$ vanishes.  This multilinear
$4$-fold product is defined as
\begin{equation*}
\imag (a \times b \times c \times d) = \imag \left( \bar a ( b ( \bar
c d) ) \right)
\end{equation*}
when $a,b,c,d$ are {\em orthogonal} octonions and $\bar a$ is the
conjugate of $a$. For non-orthogonal arguments we can write them in
terms of an orthogonal basis and expand by
multilinearity. (See~\cite{HL} Section IV.1.C for details.) Without
loss of generality we can assume that at the point $\mathbf x_0$, we
have chosen our coordinates so that $e^1 = \oee$ and $e^2 = \oie$ with
respect to the identification $T_x (\spm (\mathbb R^4)) \cong \mathbb
O$, where $\rest{T (\mathbb R^4)}{M} \cong \mathbb H \oee$ and the
spinor space $\spm \cong \mathbb H$. Similarly we can also take $\nu^1
= \oje, \nu^2 = \oke$. From this choice it follows that $\jm = \oee
(\oie) = \oi$. Then the orthonormal basis for $V_{+\jm}$ is just $q_1
= \oo, q_2 = \oi$. Now we compute (using the octonion multiplication
table):
\begin{eqnarray*}
\gamma(e^1)\gamma(\nu^1)q_1 = \oj & \quad \quad &
\gamma(e^1)\gamma(\nu^1)q_2 = \ok \\ \gamma(e^1)\gamma(\nu^2)q_1 =
\ok & \quad \quad & \gamma(e^1)\gamma(\nu^2)q_2 = -\oj \\
\gamma(\nu^1)\gamma(e^2)q_1 = \ok & \quad \quad &
\gamma(\nu^1)\gamma(e^2)q_2 = -\oj \\ \gamma(\nu^2)\gamma(e^2)q_1
= -\oj & \quad \quad & \gamma(\nu^2)\gamma(e^2)q_2 = -\ok
\end{eqnarray*}
Therefore the tangent vectors to the
immersion at
$(\mathbf x_0, t_1, t_2)$ are given by
\begin{eqnarray*}
E_1 & = & \oee + \frac{t_1}{2} \oi \left( (a_{12} -b_{11})
\oj + (a_{11} + b_{12} ) \ok \right) + \frac{t_2}{2} \oi \left(
(-a_{11} - b_{12} ) \oj + (a_{12} -b_{11}) \ok \right) \\
E_2 & = & \oie + \frac{t_1}{2} \oi \left( (a_{22} -b_{12}) \oj +
(a_{12} + b_{22} ) \ok \right) + \frac{t_2}{2} \oi \left( (-a_{12}
- b_{22}) \oj + (a_{22} - b_{12} ) \ok \right) \\ F_1 & = & \oo \\
F_2 & = & \oi
\end{eqnarray*}
Now we can compute $\imag (E_1 \times E_2 \times F_1 \times F_2)$,
with the result being
\begin{equation*}
\left( \frac{t_1}{2} (a_{11} + a_{22}) - \frac{t_2}{2} ( b_{11} +
b_{22} )\right) \oje + \left( \frac{t_1}{2} (b_{11} + b_{22}) +
\frac{t_2}{2} ( a_{11} + a_{22} )\right) \oke 
\end{equation*}
which vanishes for all $t_1, t_2$ if and only if $\tr A^{\nu_1} =
\tr A^{\nu_2} = 0$.
\end{proof}

Although this construction does produce two distinct Cayley
submanifolds of $\mathbb R^8$ for each minimal surface $M^2$ in
$\mathbb R^4$, they are in a sense degenerate examples. Note that
when the global identification of $\mathbb R^8 = \mathbb O$ has
been made, then no matter what surface $M$ we choose, the octonion
$\oo$ will be in $V_{+ \jm}$ and the space $V_{-\jm}$ will be
orthogonal to $\oo$. Therefore the $V_{+\jm}$ Cayley submanifold
will always be of the form $\mathbb R \oo \times L^3$ for some
$3$-manifold $L^3$ which therefore must be associative in $\imag
(\mathbb O) = \mathbb R^7$. Similarly the $V_{-\jm}$ Cayley
submanifold will have zero projection onto the $\oo$ component,
and thus is actually a coassociative submanifold of $\mathbb R^7$.
Note however that this does indeed give coassociative submanifolds
which are not contained in a strictly smaller subspace of $\mathbb
R^7$, which we were unable to find in Section~\ref{coasssec}. We
present some explicit examples in Section~\ref{examplessec}.
\begin{rmk}
On more general non-compact manifolds of holonomy \SP, like $\spm
(S^4)$ (see~\cite{BS, GPP}), this construction does produce
interesting Cayley submanifolds. This is discussed in~\cite{KM}.
\end{rmk}

\subsection{The space $\spi (\mathbb R^3)$ as a manifold with a
parallel \Gs} \label{spinors3sec}

A \Gs\ can similarly be placed on the spinor bundle $\spi
(\mathbb R^3) \cong \mathbb R^7$ of $\mathbb R^3$. (See~\cite{BS}
for details.) In this case we do not have positive and negative
spinor bundles. The fibre (spinor space) at each point is again
isomorphic to the quaternions $\mathbb H$. In fact we have
\begin{equation*}
\spi (\mathbb R^3) = \rest{\spi_{\pm} (\mathbb R^4)}{\mathbb R^3}
\end{equation*}
Explicitly, if $e^0, e^1, e^2, e^3$ is a basis for the Clifford
algebra of $\mathbb R^4$, then the Clifford products $e^0 \cdot
e^1$, $e^0 \cdot e^2$, $e^0 \cdot e^3$ are a basis for the Clifford
algebra of $\mathbb R^3$. We can take a surface $M^2 \subset
\mathbb R^3$ with orthonormal cotangent frame $e^1, e^2$ and
conormal vector $\nu = e^3$ and again consider the eigenspaces
$V_{\pm \jm}$ of the operator $r = \gamma(e^1) \gamma (e^2) = \pm
\gamma(e^0) \gamma(e^3)$ where the sign depends on the choice of
orientation and does not affect the eigenspaces. Then we can
take the total spaces of $V_{\pm \jm}$ over $M^2$ as $4$-manifolds
which can be coassociative in $\mathbb R^7$.
\begin{prop} \label{spinors3prop}
The total spaces of $V_{\pm \jm}$ over $M^2$ are coassociative in
$\mathbb R^7$ iff $M^2 \subset \mathbb R^3$ is {\em minimal}.
\end{prop}
\begin{proof}
Since being coassociative in $\mathbb R^7$ is equivalent to
being Cayley in $\mathbb R^8$, Theorem~\ref{cayleythm} says that
$M^2$ must be minimal in $\mathbb R^4 = \mathbb R \times \mathbb
R^3$. But since $M^2$ sits in $\mathbb R^3 \subset \mathbb R^4$,
this is equivalent to being minimal in $\mathbb R^3$.
\end{proof}

Similarly we can try to take a curve $C^1 \subset \mathbb R^3$ and
decompose the spinor space $\spi$ into eigenspaces of $r = \gamma
(e^0 ) \gamma (e^1) = \pm \gamma (\nu^1) \gamma (\nu^2)$, where
$e^1$ is a unit cotangent vector to $C^1$ and $\nu^1, \nu^2$ are an
orthonormal basis of conormal vector fields. Then the total spaces
of the bundles over $C^1$ would be $3$-manifolds which could be
associative. But since $C^1$ would have to be minimal, it is a
straight line and this construction only produces associative
$3$-planes in $\mathbb R^7$.

\section{Some Explicit Examples} \label{examplessec}

\subsection{Some Explicit Minimal Surfaces in $\mathbb R^4$}
\label{minimalssec}

For the convenience of the reader, we present some explicit
examples of minimal surfaces in $\mathbb R^4$ which are used to
construct examples of calibrated submanifolds of $\mathbb R^7$ and
$\mathbb R^8$ in Section~\ref{calibexamplessec}. If we consider a
{\em graph} of the form
\begin{equation*}
\left( x^1, x^2, f^1(x^1,x^2), f^2(x^1,x^2) \right)
\end{equation*}
then the tangent vectors to this immersion are
\begin{equation*}
e_1 = \left( 1, 0 , f^1_1, f^2_1 \right) \qquad \qquad e_2 =
\left( 0, 1, f^1_2, f^2_2 \right)
\end{equation*}
where the subscript $k$ denotes partial differentiation with
respect to $x^k$. The induced metric is $g_{ij} = e_i
\cdot e_j$. The minimal surface equations in these
coordinates are
\begin{equation} \label{minimaleq}
g_{22} f^k_{11} + g_{11} f^k_{22} - 2 g_{12} f^k_{12} = 0 \qquad k
= 1,2
\end{equation}
They are a pair of second order, quasi-linear PDE's in which the
second order derivatives are uncoupled.

Let us identify $\mathbb R^4 = \mathbb C^2$ with complex
coordinates $z = x^1 + i x^2$ and $w = f^1 + i f^2$. It is well
known (and trivial to check) that the image of a holomorphic or
anti-holomorphic map $w = f(z)$ is a minimal surface. These
satisfy the Cauchy-Riemann equations $f^1_1 = f^2_2$ and $f^1_2 =
-f^2_1$ in the holomorphic case and $f^1_1 = -f^2_2$ and $f^1_2 =
f^2_1$ in the anti-holomorphic case.

Alternatively we can instead choose complex coordinates $z=x^1 + i
f^1$ and $w = x^2 + i f^2$. Then a special Lagrangian graph is an
example of a minimal surface in $\mathbb R^4 = \mathbb C^2$. In this
case $f^k = \frac{\partial F}{\partial x^k}$ for some potential
function $F(x^1, x^2)$ and the special Lagrangian differential
equation with phase $e^{i \theta}$ is
\begin{eqnarray} \label{slagseq}
F_{11} + F_{22} & = & 0 \qquad \qquad \text{ for } \theta =
0 \\ \nonumber F_{11} F_{22} - F_{12}^2 & = & 1 \qquad \qquad
\text{ for } \theta = \frac{\pi}{2} 
\end{eqnarray}

We can also look for minimal surfaces which are not of these
special types. Our first example is a generalization of the
holomorphic example $f^1 = e^u \cos(v)$, $f^2 = e^u \sin(v)$,
which corresponds to the holomorphic function $e^z$ where we
are now writing $z = u + i v$. We can ask for the most general
minimal surface of the form
\begin{equation*}
\left( u, v, f(u) \cos(v), f(u) \sin(v) \right)
\end{equation*}
for some function $f(u)$. Substitution into~\eqref{minimaleq}
yields the following non-linear ODE for $f(u)$:
\begin{equation*}
f ( 1 + (f')^2 ) = f'' ( 1 + f^2 )
\end{equation*}
This can be explicitly integrated to give the general solution
\begin{equation*}
f(u) = \frac{C}{2}e^{K u} + \frac{1 - K^2}{2 C K^2} e^{-K u}
\end{equation*}
for two constants of integration $C$ and $K$. Note that $K=1$
corresponds to the holomorphic solution $e^u$. In
Section~\ref{calibexamplessec} we use this minimal surface with
$C=2$ and $K = \frac{1}{2}$:
\begin{equation} \label{nontriv1eq}
\left( u, v, \left (e^{\frac{u}{2}} + \frac{3}{4} e^{-\frac{u}{2}}
\right) \cos(v), \left( e^{\frac{u}{2}} + \frac{3}{4}
e^{-\frac{u}{2}} \right) \sin(v) \right)
\end{equation}

Another explicit example can be obtained by considering graphs
which are rotationally symmetric:
\begin{equation*}
\left( u, v, f(u^2 + v^2) , g(u^2 + v^2) \right)
\end{equation*}
This time substitution into~\eqref{minimaleq} yields the
following system of non-linear ODE's, where we have denoted $t =
u^2 + v^2$:
\begin{eqnarray*}
t f'' + f' + 2 t f' \left( (f')^2 + (g')^2 \right) & = & 0 \\ t g''
+ g' + 2 t g' \left( (f')^2 + (g')^2 \right) & = & 0 
\end{eqnarray*}
These can also be integrated explicitly to obtain
\begin{eqnarray*}
f(t) & = &  \frac{2 K}{\sqrt{L}} \log \left( \sqrt{t} + \sqrt{t -
\frac{4 ( 1 + K^2)}{L}} \right) \\ g(t) & = & \frac{2}{\sqrt{L}}
\log \left( \sqrt{t} + \sqrt{t - \frac{4 ( 1 + K^2)}{L}} \right)
\end{eqnarray*}
for two constants of integration $K$ and $L$. Note that this
example is only defined outside a circle in the $u,v$ plane. We use
this minimal surface in Section~\ref{calibexamplessec} with
$K=1$ and $L=4$:
\begin{equation} \label{nontriv2eq}
\left( u, v, \log \left( \sqrt{u^2 + v^2} + \sqrt{u^2 + v^2 - 2}
\right) , \log \left( \sqrt{u^2 + v^2} + \sqrt{u^2 + v^2 - 2}
\right) \right)
\end{equation}

\subsection{Examples of Calibrated Submanifolds}
\label{calibexamplessec}

We now apply the constructions described in
Section~\ref{exceptionalsec} to some explicit examples. Our
surfaces $M^2$ will all be given as graphs $(u, v, f^1(u,v),
f^2(u,v))$.

It can be checked easily that anti-holomorphic surfaces (or
equivalently special Lagrangian surfaces of any phase) satisfy the
real isotropic minimal surface equation (with the minus sign) from
Theorem~\ref{coassthm} that was required to construct coassociative
submanifolds. One can check that in these cases the constructed
$4$-fold is simply a product $\mathbb R^2 \times M^2$. Similarly a
product $3$-manifold $\mathbb R \times M^2$ is obtained when using
these minimal surfaces to construct associative submanifolds using
Theorem~\ref{assthm}.

However, we can also try holomorphic surfaces (which are still
minimal) in the associative case. (Recall that these satisfy the
real isotropic equation with the plus sign, and cannot be used to
construct coassociative submanifolds. They would work in
$\wedge^2_+ (\mathbb R^4)$, but would produce product manifolds
there.) Consider the holomorphic surface $(x, y, u(x,y), v(x,y))$
in $\mathbb R^4$ where the Cauchy-Riemann equations $u_x = v_y$
and $u_y = - v_x$ are satisfied. Then one can construct the vector
$e^1 \wedge e^2 - \nu^1 \wedge \nu^2$ in $\wedge^2_-$ and it turns
out to be (using the Cauchy-Riemann equations to simplify):
\begin{equation*}
\left( \frac{ 1 - {|\nabla u|}^2}{1 + {|\nabla u|}^2} , \frac{2
u_y}{1 + {|\nabla u|}^2}, \frac{2 u_x}{1 + {|\nabla u|}^2} \right)
\end{equation*}
Hence Theorem~\ref{assthm} gives the following associative
submanifold of $\mathbb R^7$:
\begin{equation*}
\left( t\frac{ 1 - {|\nabla u|}^2}{1 + {|\nabla u|}^2} , t\frac{2
u_y}{1 + {|\nabla u|}^2}, t\frac{2 u_x}{1 + {|\nabla u|}^2} ,
x, y, u(x,y), v(x,y) \right)
\end{equation*}
For an explicit example, we can take $u = e^x \cos(y)$ and $v =
e^x \sin(y)$ to obtain
\begin{equation*}
\left( t\frac{\sinh(x)}{\cosh(x)} , t\frac{\sin(y)}{\cosh(x)},
-t\frac{\cos(y)}{\cosh(x)} , x, y, e^x \cos(y), e^x
\sin(y) \right)
\end{equation*}
If we take instead the minimal surface in~\eqref{nontriv1eq} we
obtain, after rescaling the fibre direction basis vector to
simplify the expression, the following non-trivial associative
submanifold of $\mathbb R^7$,:
\begin{equation*}
\left( t\frac{4 e^x - 9}{12 e^{\frac{1}{2}x}} , t\sin(y), -t\cos(y)
, x, y, \left (e^{\frac{x}{2}} + \frac{3}{4} e^{-\frac{x}{2}}
\right) \cos(y), \left (e^{\frac{x}{2}} + \frac{3}{4}
e^{-\frac{x}{2}} \right) \sin(y) \right)
\end{equation*}
Finally, the minimal surface in~\eqref{nontriv2eq} yields the
following associative submanifold of $\mathbb R^7$ (defined for
$x^2 + y^2 > 2$):
\begin{equation*}
\left( (y - x)h_1 h_2, y - x, x + y, x, y, \log \left( h_1 +
h_2 \right) , \log \left( h_1 + h_2 \right) \right)
\end{equation*}
where $h_1(x,y) = \sqrt{x^2 + y^2}$ and $h_2(x,y) = \sqrt{x^2 + y^2
- 2}$.

Recall from the remarks made at the end of Section~\ref{cayleysec}
that the Cayley construction actually produces Cayley submanifolds
which are either a line cross an associative submanifold of
$\mathbb R^7$ or a coassociative submanifold of $\mathbb R^7$.
Thus they can be used to provide non-trivial examples
of coassociative submanifolds which are not contained in a
strictly smaller subspace of $\mathbb R^7$, by taking the
$V_{-\jm}$ eigenspace. Taking a holomorphic surface
$(x, y, u(x,y), v(x,y))$ in $\mathbb R^4$, one can compute that the
$-\jm$ eigenspace is spanned by
\begin{equation*}
(0, - 2 u_y, 1 - {|\nabla u|}^2, 0) \qquad \text{ and } \qquad (0,
- 2 u_x, 0, 1 - {|\nabla u|}^2)
\end{equation*}
Thus Theorem~\ref{cayleythm} gives the following coassociative
submanifold of $\mathbb R^7$:
\begin{equation*}
\left( -2 (t_1 u_y + t_2 u_x) , t_1(1 - {|\nabla u|}^2), t_2(1 -
{|\nabla u|}^2) , x, y, u(x,y), v(x,y) \right)
\end{equation*}
The example of $u = e^x \cos(y)$ and $v = e^x \sin(y)$ gives
\begin{equation*}
\left( 2 e^x (t_1 \sin(y) - t_2 \cos(y)) , t_1 (1 - e^{2 x}), t_2 (1 -
e^{2 x}) , x, y, e^x \cos(y), e^x \sin(y) \right)
\end{equation*}
as a coassociative submanifold of $\mathbb R^7$. One can similarly
use~\eqref{nontriv1eq} or~\eqref{nontriv2eq} and
Theorem~\ref{cayleythm} to produce explicit coassociative
submanifolds of $\mathbb R^7$. The expressions tend to be
extremely complicated in these cases.

\appendix

\section{Octonion Multiplication Table} \label{octtablesec}

The following is a multiplication table for the octonions $\mathbb
O$. The table corresponds to multiplying the element in the
corresponding row on the left of the element in the corresponding
column. For example $\oi \cdot \oj = \ok$.

\begin{table}[h]
\begin{center}
\begin{tabular}{|c||c|c|c|c|c|c|c|c|} 
\hline {} & \oo & \oi & \oj & \ok & \oee & \oie & \oje & \oke \\
\hline \hline \oo & \oo & \oi & \oj & \ok & \oee & \oie & \oje &
\oke \\ \hline \oi & \oi & -\oo & \ok & -\oj & \oie & -\oee & -\oke
& \oje \\ \hline \oj & \oj & -\ok & -\oo & \oi & \oje & \oke &
-\oee & -\oie \\ \hline \ok & \ok & \oj & -\oi & -\oo & \oke &
-\oje & \oie & -\oee \\ \hline \oee & \oee & -\oie & -\oje & -\oke
& -\oo & \oi & \oj & \ok \\ \hline \oie & \oie & \oee & -\oke &
\oje & -\oi & -\oo & -\ok & \oj \\ \hline \oje & \oje & \oke &
\oee & -\oie & -\oj & \ok & -\oo & - \oi \\ \hline \oke & \oke &
-\oje & \oie & \oee & -\ok & -\oj & \oi & -\oo \\ \hline
\end{tabular}
\end{center}
\end{table}

\end{document}